\documentclass[11pt,twoside]{article}
\usepackage{latexsym,amsthm,amsmath,amsfonts,amscd,amssymb}

\evensidemargin\oddsidemargin
\theoremstyle{plain}  % default
\newtheorem{theorem}{Theorem}[section]

\newtheorem*{theorem*}{Theorem}

\newtheorem{lemma}[theorem]{Lemma}
\newtheorem{proposition}[theorem]{Proposition}

\newtheorem{tech-lemma}[theorem]{Technical Lemma}
\newtheorem{definition}[theorem]{Definition}
\theoremstyle{remark}

\newtheorem{remark}[theorem]{Remark}
\newtheorem*{remark*}{Remark}

\newtheorem*{claim*}{Claim}
\numberwithin{equation}{section}
\renewcommand{\leq}{\leqslant}

\renewcommand{\geq}{\geqslant}

\newcommand{\lto}{\longrightarrow}
\newcommand{\lmapsto}{\longmapsto}
\newcommand{\R}{\mathbb{R}}
\newcommand{\Z}{\mathbb{Z}}
\newcommand{\C}{\mathbb{C}}

\newcommand{\norm}[1]{\lVert#1\rVert}
\newcommand{\lie}{\mathfrak}

\newcommand{\suchthat}{\;\;|\;\;}

\newcommand{\PGL}{\mathrm{PGL}} 
\newcommand{\PSL}{\mathrm{PSL}}

\newcommand{\GL}{\mathrm{GL}}

\newcommand{\g}{\mathfrak{g}}
\newcommand{\h}{\mathfrak{h}}
\newcommand{\m}{\mathfrak{m}}
\newcommand{\SL}{\mathrm{SL}}

\newcommand{\Or}{\mathrm{O}}

\newcommand{\la}{\langle}
\newcommand{\ra}{\rangle}

\DeclareMathOperator{\ad}{ad}

\DeclareMathOperator{\rk}{rk}

\DeclareMathOperator{\Hom}{Hom}
\DeclareMathOperator{\End}{End}

\begin{document}
\null\vskip 10mm

%%%%%%%%%%%%%%%%%%%%%%%%%%%%%%%%%%%%%%%%%%%%%%%%%%%%%%%%%%%%%%%%
%
% Data for the headings. Please fill in both fields.
%
%%%%%%%%%%%%%%%%%%%%%%%%%%%%%%%%%%%%%%%%%%%%%%%%%%%%%%%%%%%%%%%%

\pagestyle{myheadings}
\markboth
{{\small\sc Representations of surface groups}}
{{\small\sc Bradlow, Garc{\'\i}a--Prada \& Gothen}}

%%%%%%%%%%%%%%%%%%%%%%%%%%%%%%%%%%%%%%%%%%%%%%%%%%%%%%%%%%%%%%%%
%
% Title. For instance:
%
%%%%%%%%%%%%%%%%%%%%%%%%%%%%%%%%%%%%%%%%%%%%%%%%%%%%%%%%%%%%%%%%

\thispagestyle{empty}

\begin{center}
{\Large\bf
Representations of surface groups 
\medskip

in the general linear group
}
\vskip 10mm

%%%%%%%%%%%%%%%%%%%%%%%%%%%%%%%%%%%%%%%%%%%%%%%%%%%%%%%%%%%%%%%%
%
% Author(s), affiliation(s) and email(s). For instance:
%
%%%%%%%%%%%%%%%%%%%%%%%%%%%%%%%%%%%%%%%%%%%%%%%%%%%%%%%%%%%%%%%%

{\large\bf
$^1$Steven B. Bradlow, $^2$Oscar Garc{\'\i}a--Prada and \\ $^3$Peter
B. Gothen
}
\vskip 5mm
{\it
$^1$Department of Mathematics,
University of Illinois,
Urbana,
IL 61801, 
USA 
\smallskip

$^2$Instituto de Matem\'aticas y F\'{\i}sica Fundamental,
Consejo Superior de Investigaciones Cient\'{\i}ficas,
Serrano 113 bis,
28006 Madrid, Spain

\smallskip

$^3$
Departamento de Matem{\'a}tica Pura,
Faculdade de Ci{\^e}ncias,
Universidade do Porto,
Rua do Campo Alegre 687, 4169-007 Porto,
Portugal

\smallskip

emails: {bradlow@math.uiuc.edu, oscar.garcia-prada@uam.es,
  pbgothen@fc.up.pt}
}
\end{center}

\bigskip

%%%%%%%%%%%%%%%%%%%%%%%%%%%%%%%%%%%%%%%%%%%%%%%%%%%%%%%%%%%%%%%%
%
% Abstract, key words and MSC codes
%
%%%%%%%%%%%%%%%%%%%%%%%%%%%%%%%%%%%%%%%%%%%%%%%%%%%%%%%%%%%%%%%%

\begin{abstract}
\parindent0pt\noindent
We determine the number of connected components of the moduli space 
for representations of a surface group in the general linear group. 



\end{abstract}

%%%%%%%%%%%%%%%%%%%%%%%%%%%%%%%%%%%%%%%%%%%%%%%%%%%%%%%%%%%%%%%%
%
% Text.
%
%%%%%%%%%%%%%%%%%%%%%%%%%%%%%%%%%%%%%%%%%%%%%%%%%%%%%%%%%%%%%%%%

%%%%%%%%%%%%%%%%%%%%%%%%%%%%%%%%%%%%%%%%%%%%%%%%%%%%%%%%%%%%%%%%%%%%%%
\section{Introduction}
%%%%%%%%%%%%%%%%%%%%%%%%%%%%%%%%%%%%%%%%%%%%%%%%%%%%%%%%%%%%%%%%%%%%%%

Let $X$ be a closed oriented surface of genus $g \geq 2$ and let
$\mathcal{M}^{+}_{G}$ be the moduli space of semisimple
representations of $\pi_{1}(X)$ in $\GL(n,\R)$.  In this paper we 
extend the results of Hitchin 
\cite{hitchin:1992} on the number of connected components of 
$\mathcal{M}^{+}_{\PSL(n,\R)}$ to the case of
$\mathcal{M}^{+}_{\GL(n,\R)}$ when $n \geq 3$. 
(The case of representations in $\PSL(2,\R)$ was studied by Goldman
\cite{goldman:1988} and the case of representations in $\PGL(2,\R)$
and $\GL(2,\R)$ was studied by Xia \cite{xia:1997, xia:1999}.)

We adopt the Morse theoretic approach pioneered by Hitchin in
\cite{hitchin:1987,hitchin:1992} and, indeed, most of our arguments
follow \cite{hitchin:1992} quite closely.  However, at one point
(Theorem~\ref{thm:GL-minima} below), the application of a general
result of \cite{BGG:2003a} allows for a significant simplification of
the arguments (cf.\ \cite[Lemma~9.6]{hitchin:1992}). This provides a
nice example of the power of the machinery introduced in that paper.

It is a pleasure to thank Nigel Hitchin for enlightening conversations
on this topic.  We also thank the referee for useful comments.

%%%%%%%%%%%%%%%%%%%%%%%%%%%%%%%%%%%%%%%%%%%%%%%%%%%%%%%%%%%%%%%%%%%%%%
\section{Representations and the moduli space}
%%%%%%%%%%%%%%%%%%%%%%%%%%%%%%%%%%%%%%%%%%%%%%%%%%%%%%%%%%%%%%%%%%%%%%

Let $X$ be a closed oriented surface of genus $g$ and let 
\begin{displaymath}
  \pi_{1}(X) = \la a_{1},b_{1}, \dotsc, a_{g},b_{g} \suchthat 
  \prod_{i=1}^{g}[a_{i},b_{i}] = 1 \ra
\end{displaymath}
be its fundamental group.  By a \emph{representation} of $\pi_1(X)$ in
$\GL(n,\R)$ we understand a homomorphism $\rho\colon \pi_1(X) \to
\GL(n,\R)$.  The set of all such homomorphisms,
$\Hom(\pi_1(X),\GL(n,\R))$, can be naturally identified with the subset
of $\GL(n,\R)^{k}$ consisting of $2g$-tuples
$(A_{1},B_{1}\dotsc,A_{g},B_{g})$ satisfying the algebraic equation
$\prod_{i=1}^{g}[A_{i},B_{i}] = 1$.  This shows that
$\Hom(\pi_1(X),\GL(n,\R))$ is a real algebraic variety.

The group $\GL(n,\R)$ acts on $\Hom(\pi_1(X),\GL(n,\R))$ by conjugation:
\[
(g \cdot \rho)(\gamma) = g \rho(\gamma) g^{-1}
\]
for $g \in \GL(n,\R)$, $\rho \in \Hom(\pi_1(X),\GL(n,\R))$ and
$\gamma\in \pi_1(X)$. If we restrict the action to the subspace
$\Hom^+(\pi_1(X),\GL(n,\R))$ consisting of semi-simple
representations, the orbit space is Hausdorff.  Define the
\emph{moduli space} for representations of $\pi_1(X)$ in $\GL(n,\R)$
to be the orbit space
\[
\mathcal{M}^{+}_{\GL(n,\R)} = \Hom^{+}(\pi_1(X),\GL(n,\R)) / \GL(n,\R)
\]
with the quotient topology.

Given a representation $\rho\colon\pi_{1}(X) \to
\GL(n,\R)$, there is an associated flat bundle on
$X$, defined as
\begin{math}
  V_{\rho} = \tilde{X}\times_{\pi_{1}(X)}\R^n
\end{math},
where $\tilde{X} \to X$ is the universal cover and $\pi_{1}(X)$ acts 
on $\R^{n}$ via $\rho$.  We then define invariants of $\rho$ as the 
Stiefel--Whitney classes of $V_{\rho}$:
\begin{align*}
  w_{1}(\rho) &= w_{1}(V_{\rho}) \in H^{1}(X,\Z/2)\ , \\
  w_{2}(\rho) &= w_{2}(V_{\rho}) \in H^{2}(X,\Z/2)\ .
\end{align*}
For fixed $(w_{1},w_{2}) \in H^{1}(X,\Z/2) \oplus H^{2}(X,\Z/2)$ we 
define a subspace
%have a corresponding subspace of $\mathcal{M}^{+}_{\GL(n,\R)}$,
\begin{displaymath}
  \mathcal{M}^{+}_{\GL(n,\R)}(w_{1},w_{2}) = \{\rho \suchthat 
  w_{i}(\rho) = w_{i},\; i=1,2 \}
  \subseteq \mathcal{M}^{+}_{\GL(n,\R)}\ .
\end{displaymath}

%%%%%%%%%%%%%%%%%%%%%%%%%%%%%%%%%%%%%%%%%%%%%%%%%%%%%%%%%%%%%%%%%%%%%%
\section{$\GL(n,\R)$-Higgs bundles}
%%%%%%%%%%%%%%%%%%%%%%%%%%%%%%%%%%%%%%%%%%%%%%%%%%%%%%%%%%%%%%%%%%%%%%

Let $G$ be a real reductive Lie group, let $H \subset G$ be a
maximal compact subgroup and let $G_{\C}$ be the complexification of
$G$.  The complexification of $H$ is denoted by $H_{\C} \subset
G_{\C}$.  At the Lie algebra level we have the Cartan decomposition
$\lie{g} = \lie{h} \oplus \lie{m}$.  The restriction of the adjoint
representation to $H$ defines a representation on $\lie{m}$ called the
\emph{isotropy representation}; let
\begin{equation}\label{eq:isotropy}
  \iota\colon H_{\C}\longrightarrow \GL(\lie{m}_{\C})
\end{equation}
be its complexification.  There is a
complex linear Lie algebra involution $\theta\colon \lie{g}_{\C} \to
\lie{g}_{\C}$ which has $\lie{h}_{\C}$ as its $+1$-eigenspace and
$\lie{m}_{\C}$ as its $-1$-eigenspace, giving the decomposition
$\lie{g}_{\C} = \lie{h}_{\C} \oplus \lie{m}_{\C}$.

In the case of interest to us, namely that of $G=\GL(n,\R)$, the
maximal compact subgroup is $H = \Or(n)$ and its complexification is
$H_{\C} = \Or(n,\C)$.  In terms of a defining representation, we thus
have a non-degenerate quadratic form $Q$ on $\C^n$, and $\Or(n,\C)$ is
the group of complex linear automorphisms of $\C^n$ preserving the
form $Q$.  The corresponding decomposition at the Lie algebra level is
$\lie{gl}(n,\C) = \lie{o}(n,\C) \oplus \m_{\C}$, where $\lie{o}(n,\C)$
and $\m_{\C}$ are respectively the antisymmetric and symmetric
endomorphisms of $\C^n$ (with respect to $Q$).  

Fix a complex structure on $X$.  Since no confusion is likely to
arise, we shall also call the corresponding Riemann surface 
$X$.  We denote the canonical bundle of $X$ by $K = T^*X^{1,0}$. 
One more piece of convenient notation is the following: for any Lie
group $G$, if we have a principal $G$-bundle $E$ on $X$ and a
$G$-space $V$, we denote the associated bundle with fibres $V$ by
$E(V) = E \times_{G} V$.  For example, in the case of the representation 
\eqref{eq:isotropy} we obtain a bundle $E(\lie{m}_{\C}) = E 
\times_{\iota}\lie{m}_{\C}$ with fibres isomorphic to $\lie{m}_{\C}$.
\begin{definition}
  A \emph{$G$-Higgs bundle} is a pair $(E,\Phi)$, where $E \to X$ is a 
  principal holomorphic $H_{\C}$-bundle and the \emph{Higgs field} 
  $\Phi$ belongs to $H^0(E(\lie{m}_{\C})\otimes K)$.
\end{definition}

A $\GL(n,\R)$-Higgs bundle is thus a pair $(E,\Phi)$, where $E$ is a
$\Or(n,\C)$ principal bundle and $\Phi \in H^0(E(\m_{\C})\otimes K)$.
Letting $V = E \times_{\Or(n,\C)}\C^n$ be the associated vector
bundle, we can describe this more concretely as a triple $(V,Q,\Phi)$,
where $V$ is a holomorphic rank $n$ vector bundle, $Q \in H^0(S^2V^*)$
is a non-degenerate quadratic form and $\Phi\in H^0(\End(V)\otimes K)$
is symmetric with respect to $Q$.  We denote by $q\colon V \to V^*$
the symmetric isomorphism associated to $Q$.
In the decomposition
\begin{displaymath}%\label{eq:decomp-EndV}
  E(\g_{\C}) = E(\h_{\C}) \oplus E(\m_{\C})
\end{displaymath}
we have that $E(\g_{\C})$ is just $\End(V)$, and the involution
$\theta$ on $\End(V)$ defining the decomposition is
\begin{equation}\label{eq:theta}
  \theta\colon A \lmapsto -(qAq^{-1})^t\ .  
\end{equation}
Thus, under the isomorphism $q\colon V \to V^*$, we can identify the
$+1$-eigenbundle of $\theta$,
\begin{displaymath}
  U^+ = E(\h_{\C})\ ,
\end{displaymath}
as the bundle of antisymmetric endomorphisms of $V$
and the $-1$-eigenbundle of $\theta$,
\begin{displaymath}
  U^- = E(\m_{\C})\ ,
\end{displaymath}
as the bundle of symmetric endomorphisms of $V$.

The notion of isomorphism between $G$-Higgs bundles is the obvious
one: $(E,\Phi)$ and $(E',\Phi')$ are isomorphic, if there is an
isomorphism $g \colon E \overset{\cong}{\lto} E'$ which takes $\Phi$
to $\Phi'$ under the induced isomorphism $E(\m_{\C})
\overset{\cong}{\lto} E'(\m_{\C})$.

There is a stability condition for $G$-Higgs bundles, which
generalizes the usual stability condition for Higgs vector bundles
(see \cite{hitchin:1987,simpson:1992}) and Ra\-ma\-na\-than's stability
condition for principal bundles \cite{ramanathan}. It is a special case
of the very general stability notion studied in \cite{BGM}.  We shall
not need the detailed form of this stability condition here, but we
mention that in the case of $\GL(n,\R)$-Higgs bundles $(V,Q,\Phi)$, it
means that the usual slope condition is satisfied for $\Phi$-invariant
isotropic subbundles of $V$.  We also note that there are
corresponding notions of poly-stability and semi-stability.

One expects this stability condition to be appropriate for
constructing moduli spaces of $G$-Higgs bundles.  For this one should
first fix the topological type of the principal $H_{\C}$-bundle $E$.
In the case of $\GL(n,\R)$-Higgs bundles, this means fixing the first
and second Stiefel--Whitney classes of any $\Or(n)$-bundle obtained 
by reduction of the structure group to the maximal compact subgroup 
$\Or(n) 
\subseteq \Or(n,\C)$.  In general no direct construction of the moduli
spaces is currently available---in the case of present interest, we
shall adopt the solution of \cite{hitchin:1992}, realizing the moduli
space of semistable $\GL(n,\R)$-Higgs bundles with fixed
Stiefel--Whitney classes $w_1$ and $w_2$ as a subspace of the moduli
space of ordinary Higgs (vector) bundles of rank $n$ and degree $0$.
We denote the moduli space by
\begin{math}
  \mathcal{M}(w_1,w_2).
\end{math}

The $G$-Higgs bundle stability condition is equivalent to an 
existence criterion for solutions to certain gauge theoretic 
equations (Hitchin's equations) on $(E,\Phi)$.  This provides a 
bridge from poly-stable $G$-Higgs bundles to representations of 
$\pi_1(X)$, since solutions to the equations give rise to flat
$G$-bundles.  The existence of solutions was proved by
Hitchin \cite{hitchin:1987} and, more generally, Simpson
\cite{simpson:1988, simpson:1992} in the case of Higgs bundles with
complex structure group.  However, in the present generality, this
requires the results of \cite{BGM}.

The correspondence in the other direction is given by the theorem of
Do\-nald\-son \cite{donaldson} and, more generally, Corlette
\cite{corlette}: given a flat semi-simple $G$-bundle, there is a
preferred reduction of structure group to the maximal compact $H
\subset G$, a so-called harmonic metric.  This gives rise to a
solution to Hitchin's equations and thus to a (poly-)stable $G$-Higgs
bundle $(E,\Phi)$; the $H_{\C}$-bundle $E$ of course being the
complexification of the $H$-bundle obtained via the harmonic metric.

We thus have the following fundamental result, essentially due to
Corlette, Donaldson, Hitchin and Simpson.

\begin{theorem}\label{thm:CDHS}
  Let $w_1 \in H^1(X,\Z/2)$ and $w_2 \in H^2(X,\Z/2)$ be fixed.  Then
  there is a homeomorphism between the moduli space of semi\-stable
  $G$-Higgs bundles $\mathcal{M}(w_1,w_2)$ and the moduli space
  $\mathcal{M}^{+}_{\GL(n,\R)}(w_{1},w_{2})$ of semisimple
  representations of the fundamental group of $X$ in $\GL(n,\R)$, with
  the given invariants.
\end{theorem}

%%%%%%%%%%%%%%%%%%%%%%%%%%%%%%%%%%%%%%%%%%%%%%%%%%%%%%%%%%%%%%%%%%%%%%
\section{Morse theory on the Higgs bundle moduli space}
%%%%%%%%%%%%%%%%%%%%%%%%%%%%%%%%%%%%%%%%%%%%%%%%%%%%%%%%%%%%%%%%%%%%%%

Theorem~\ref{thm:CDHS} shows that counting the number of connected
components of the moduli space of representations is the same thing as
counting the number of connected components of the moduli space of
$\GL(n,\R)$-Higgs bundles.  In order to do this we use the Morse
theory approach to the study of the topology of Higgs bundle moduli
introduced by Hitchin in \cite{hitchin:1987,hitchin:1992}.  The
$L^2$-norm of $\Phi$ defines a positive function
\begin{math}
  f\colon \mathcal{M}(w_1,w_2) \lto \R \end{math},
given by $f(E,\Phi) = \int_{X}\norm{\Phi}^2$ (this definition uses the
harmonic metric in the bundle $E$).  When the Higgs bundle moduli
space is smooth, this function is a perfect Bott--Morse function,
giving a powerful tool for the study of the topology of the moduli
space.  But even when singularities are present, the fact that $f$ is
a proper map gives the following result on connected components.
\begin{proposition}\label{prop:proper}
  Let $\mathcal{M}' \subseteq
  \mathcal{M}(w_1,w_2)$ be a subspace and let
  $\mathcal{N} \subseteq \mathcal{M}'$ be the subspace of
  $\mathcal{M}'$ consisting of local minima of the restriction of $f$. 
  If $\mathcal{N}$ is connected, then so is $\mathcal{M}'$.
\end{proposition}
Using this result for determining the connected components of
$\mathcal{M}(w_1,w_2)$ obviously requires
identifying the local minima of $f$ on
$\mathcal{M}(w_1,w_2)$.  In order to do this, we
use one of the main theorems proved in \cite{BGG:2003a}.  Before
stating the result we need some preliminaries.

Any local minimum of $f$ corresponds to a fixed point of the
$\C^*$-action $(E,\Phi) \mapsto (E,\lambda\Phi)$.  It is not hard to
see that $(E,\Phi)$ represents a fixed point if and only if it is a
so-called \emph{complex variation of Hodge structure}.  In the case of
a $\GL(n,\R)$-Higgs bundle $(V,Q,\Phi)$, this means the following
(cf.\ Hitchin \cite{hitchin:1992} p.\ 466): the vector bundle $V$
breaks up as a direct sum $V = F_{-m} \oplus \cdots \oplus F_{m}$ and
the restriction $\Phi_i$ of the Higgs field $\Phi$ to $F_i$, gives
maps $\Phi_i \colon F_{k} \to F_{k+1} \otimes K$.  Furthermore, the
quadratic form $Q$ gives an isomorphism 
$q\colon F_{k} \to F_{-k}^*$,
and the remaining $F_{l}$ are orthogonal to $F_{k}$ under $Q$. 
Recalling that $\Phi$ is symmetric with respect to $Q$ we thus have
that
\begin{equation}
  \label{eq:phi_k-sym}
  \Phi_{-k} = \Phi_{k-1}^t \colon F_{-k} \to F_{-k+1}\otimes K\ .
\end{equation}
Define $U_{ij} = \Hom(F_{j},F_{i})$ and $U_{k} = \bigoplus_{i-j=k}
U_{ij}$.  Then there is a corresponding decomposition of the Lie
algebra bundle $E(\lie{g}_{\C}) = \End(V)$ as
\begin{displaymath}
  \End(V) = \bigoplus_{k=-2m}^{2m} U_{k}\ .
\end{displaymath}
The restriction of the
involution $\theta$ gives an isomorphism 
\begin{math}
  \theta\colon U_{ij} \lto U_{-j,-i}
\end{math}
(cf.\ \eqref{eq:theta}).
Hence $\theta$ restricts to $\theta \colon U_{k} \to U_{k}$.  Letting $U_{k} =
U_{k}^{+} \oplus U_{k}^{-}$ be the corresponding eigenspace
decomposition, we thus have $U^{+} = \bigoplus U_{k}^{+}$ and $U^{-} =
\bigoplus U_{k}^{-}$.  The fact that $\Phi$ maps $F_{k}$ to $F_{k+1}
\otimes K$ means that $\Phi \in H^0(U_{1}^{-} \otimes K)$.  Note that
$\ad(\Phi)$ interchanges $U^{+}$ and $U^{-}$, and therefore $\ad(\Phi)
\colon U_{k}^{\pm} \to U_{k+1}^{\mp} \otimes K$.
The result we need from \cite{BGG:2003a} can now be stated as follows.
\begin{theorem}
\label{thm:adjoint-minima}
  Let $(V,Q,\Phi)$ be a stable $\GL(n,\R)$-Higgs bundle which represents a
  critical point of $f$.  Then this critical point is a local minimum 
  if and only if either $\Phi=0$ or
  \begin{displaymath}
    \ad(\Phi) \colon U_{k}^{+} \to U_{k+1}^{-} \otimes K
  \end{displaymath}
  is an isomorphism for all $k \geq 1$.
\end{theorem}

\begin{proof}
  This follows from \cite[Proposition~4.14]{BGG:2003a} by an argument 
  analogous to the proof of \cite[Corollary~4.15]{BGG:2003a} (cf.\ 
  \cite[Remark~4.16]{BGG:2003a}).
\end{proof}  

\begin{theorem}\label{thm:GL-minima}
  Let the stable $\GL(n,\R)$-Higgs bundle $(V,Q,\Phi)$ be a complex
  variation of Hodge structure.  Assume that $n \geq 3$.  Then
  $(V,Q,\Phi)$ represents a minimum of $f$ if and only if one of the
  following two alternatives occurs:
  \begin{itemize}
    \item[$(1)$] The Higgs field $\Phi$ vanishes identically.
  
    \item[$(2)$] Each bundle $F_{i}$ has rank $1$ and, furthermore,
    the restriction $\Phi_{i} = \Phi_{|F_{i}}$ defines an isomorphism
    $\Phi_{i}\colon F_{i} \overset{\cong}{\lto} F_{i+1}\otimes K$.
  \end{itemize}
\end{theorem}

\begin{remark}
  In the case $n=1$, it is easy to see that only minima of the type 
  described in $(1)$ of the Theorem occur.  The case $n=2$ is also 
  special, since a third type of minima exists (see 
  \cite[Proposition~9.19]{hitchin:1992}).
\end{remark}

\begin{proof}[{Proof of Theorem~\ref{thm:GL-minima}}]
  It is clear from the definition of $f$ as the $L^2$-norm of $\Phi$
  that $\Phi = 0$ implies that $(E,\Phi)$ is a local minimum.  So
  assume from now on that $\Phi \neq 0$ and represents a local minimum
  of $f$.
  
  Consider first $U_{m,-m} = \Hom(F_{-m},F_m) = U_{2m}$.  Then we have
  the isomorphism $\theta\colon U_{m,-m} \overset{\cong}{\lto}
  U_{m,-m}$ and the $+1$-eigenspace $U^{+}_{m,-m}$ is the space of
  antisymmetric maps $F_{-m} \to F_m$ under the duality $q\colon F_{m}
  \overset{\cong}{\lto} F_{-m}^{*}$.  Theorem~\ref{thm:adjoint-minima}
  says that we have an isomorphism $\ad(\Phi)\colon U_{2m}^{+} \to
  U_{2m+1}^{-} \otimes K$.  Since the latter space is $0$, so is the
  former.  Hence we conclude that there are no antisymmetric maps
  $F_{-m} \to F_m$.  This is only possible if $\rk(F_{m}) =
  \rk(F_{-m}) = 1$.
  
  Next we shall prove that the remaining $F_{i}$ are line bundles and
  that $\Phi_i\colon F_{i} \to F_{i+1}\otimes K$ are isomorphisms.  Note
  that we only need to do this for $i \geq 0$ (cf.\
  \eqref{eq:phi_k-sym}).  We proceed by induction, taking as induction
  hypothesis that $F_{m-l}$ has rank $1$ for $0\leq l \leq k$, and 
  show for $k \leq m-1$ that $F_{m-k-1}$ has rank $1$ and that we 
  have an isomorphism $\Phi_{m-k-1} \colon 
  F_{m-k-1}\overset{\cong}{\lto} F_{m-k} \otimes K$. 
  
  Consider $U_{m-(k+1),-m}$ and $U_{m,k+1-m}$.  These spaces are
  transformed into each other by $\theta$ and, if we organize the
  $U_{ij}$ into a matrix, they are located at opposite extremes of the
  diagonal whose elements make up $U_{2m-(k+1)}$.  Thus elements of
  the form $(a,\theta(a))$ in the direct sum $U_{m-(k+1),-m} \oplus
  U_{m,k+1-m}$ belong to $U^{+}_{2m-(k+1)}$.  By
  Theorem~\ref{thm:adjoint-minima} we know that the restriction of
  $\ad(\Phi)$ to the subspace of $U^{+}_{2m-(k+1)}$ consisting of such
  pairs $(a,\theta(a))$ is injective when $k \leq 2m-2$.  This
  condition is satisfied, since $k \leq m-1$ (this is where we need $n
  \geq 3$, in order to have $m-1 \geq 0$).  We have\footnote{In this
    formula $\theta(a)$ should be twisted by the identity on $K$.
    Since no confusion can arise, we shall discard such twisting from
    the notation.}
  \begin{equation}\label{eq:ad-phi}
    \ad(\Phi)(a,\theta(a)) = 
      \Phi_{m-(k+1)} \circ a - \theta(a) \circ \Phi_{k-m}\ . 
  \end{equation}

  If $k=0$ these two summands both lie in $U_{m,-m}$.  We calculate 
  the second one, using that $\Phi_{-m} = \theta(\Phi_{m-1}) = 
  (q\Phi_{m-1}q^{-1})^{t}$:
  \begin{displaymath}
    \theta(a) \circ \Phi_{-m} 
      = -(qaq^{-1})^{t} \circ (q\Phi_{m-1}q^{-1})^{t} 
      = -(q\Phi_{m-1}aq^{-1})^{t} 
      = \theta(a\Phi_{m-1})\ .
  \end{displaymath}
  But we have already proved that $U_{m,-m}$ is a line bundle on 
  which $\theta$ is $+1$.  Hence \eqref{eq:ad-phi} shows that when $k=0$
  \begin{displaymath}
    \ad(\Phi)(a,\theta(a)) = 
      2\Phi_{m-1} \circ a\ .
  \end{displaymath}
  Since $\ad(\Phi)$ is injective, this proves that $\Phi_{m-1}\colon
  F_{m-1} \to F_{m} \otimes K$ injects.  But $F_{m}$ is a line bundle,
  so $F_{m-1}$ must also be a line bundle and $\Phi_{m-1}$ an
  isomorphism.

  A similar argument works when $k \geq 1$; however this case is
  easier, because the summands in \eqref{eq:ad-phi} lie in different
  $U_{i,j}$ and hence we can appeal to injectivity of
  $(a,\theta(a))\mapsto \Phi_{m-k+1} \circ a$.
\end{proof}

Let $(V,Q,\Phi)$ be a local minimum of $f$ of the kind described in
$(2)$ of Theorem~\ref{thm:GL-minima}.  Using $F_{k} \cong F_{-k}^*$ we
see that $n = 2m+1$ and that $m$ is half integer when $n$ is even,
while $m$ is integer when $n$ is odd.  Hence
\begin{align*}
  V &= F_{-m} \oplus \cdots \oplus F_{-1/2} \oplus F_{1/2}
    \oplus \cdots \oplus F_{m} &\text{if $n$ is even,}\\
  V &= F_{-m} \oplus \cdots \oplus F_{0}
    \oplus \cdots \oplus F_{m} &\text{if $n$ is odd.}
\end{align*}
Thus Theorem~\ref{thm:GL-minima} leads to the following more precise
characterization of the $(V,Q,\Phi)$ with $\Phi \neq 0$ representing a
local minimum of $f$.

\begin{proposition}\label{prop:char-minima}
  Let $(V,Q,\Phi)$ be a local minimum of $f$ of the kind described in
  $(2)$ of Theorem~\ref{thm:GL-minima}.  Then the following holds.
  \begin{itemize}
  \item[$(1)$] If $n$ is even, then $F_{-1/2}^2 = K$ and the 
  remaining $F_{i}$ are uniquely determined by the choice of this 
  square root of $K$ as $F_{-1/2+k} \cong F_{-1/2} \otimes K^{-k}$.
  \item[$(2)$] If $n$ is odd, then $F_0^2 = \mathcal{O}$ and 
  $F_{k} \cong F_0 \otimes K^{-k}$ for $k \neq 0$.
  \item[$(3)$] In both cases, $(V,Q,\Phi)$ is isomorphic to a 
  $\GL(n,\R)$-Higgs bundle, where 
  \begin{displaymath}
    q = 
    \begin{pmatrix}
        0 &  \cdots& \cdots & 0 & 1 \\
        \vdots & & & \cdots & 0 \\
        \vdots & & 1 & & \vdots \\
        0 & \cdots & & & \vdots \\
        1 & 0 & \cdots& \cdots & 0
    \end{pmatrix} \qquad\text{and}\qquad
    \Phi =
    \begin{pmatrix}
      0 &\cdots &\cdots &\cdots & 0 \\
      1 & 0 & \cdots & \cdots& 0\\
      0 & 1 & 0 & \cdots & 0 \\
      \vdots  &  & \ddots & & \vdots \\
      0 & \cdots & 0 & 1 & 0
    \end{pmatrix}\ , 
  \end{displaymath}
  with respect to the decomposition $V=F_{-m} \oplus \cdots \oplus
  F_{m}$.  (Here each $F_i$ is a line bundle and $n = 2m+1$.  In the
  notation for $q\colon V \overset{\cong}{\lto} V^{*}$ we use $F_{i}
  \cong F_{-i}^{*}$ and in the notation for $\Phi$ we use $1$ for the
  canonical identification $F_0\otimes K^{-k} \cong F_0\otimes
  K^{-(k+1)}\otimes K$.)
  \end{itemize}
\end{proposition}

\begin{proof}
  $(1)$ and $(2)$ are clear, using the isomorphism $\Phi_{i}\colon
  F_{i} \overset{\cong}{\lto} F_{i+1} \otimes K$.  For $(3)$ we note
  that $(V,Q,\Phi)$ is of the form given, except that the $1$'s
  appearing in $\Phi$ are arbitrary non-zero complex scalars
  $\lambda_{1},\dotsc,\lambda_{n-1}$.  It is an easy exercise to show
  that such a $(V,Q,\Phi)$ is isomorphic to a $\GL(n,\R)$-Higgs bundle
  of the kind given.
\end{proof}

Now we calculate the invariants $(w_1,w_2)$ of the $(V,Q,\Phi)$
described in the preceding Proposition.  To be precise, we need to
calculate the Stiefel--Whitney classes of the real bundle obtained by
a reduction of structure group in $(V,Q)$ from $\Or(n,\C)$ to
$\Or(n)$.  We shall denote these classes by $w_{i}(V,Q)$.

\goodbreak

\begin{proposition}\label{prop:inv-minima}
  Let $(V,Q,\Phi)$ be a local minimum of $f$ of the kind described in
  $(2)$ of Theorem~\ref{thm:GL-minima}.  If $n=2q$ is even, then
  \begin{align*}
    w_1(V,Q) &= 0 \ ,\\
    w_2(V,Q) &= (g-1)q^2 \mod 2\ ,
  \intertext{If $n=2q+1$ is odd, then}
    w_1(V,Q) &= w_1(F_0,Q_{|F_0}) \ ,\\
    w_2(V,Q) &= 0\ .
  \end{align*}
\end{proposition}

\begin{proof}
In the case of the first Stiefel-Whitney class $w_{1}$, consider the
determinant bundle $\Lambda^{n}V$ in the Jacobian of $X$.  Since $V
\cong V^{*}$, we have $(\Lambda^{n}V)^{2} = \mathcal{O}$.  It is then
easy to see that $\Lambda^{n}V$ corresponds to $w_{1}(V)$ under the
identification of the $2$-torsion points in the Jacobian with
$H^1(X,\Z/2)$.  From this the calculation of $w_{1}(V,Q)$ is 
immediate.  

For the calculation of $w_{2}$ we can get a reduction of structure
group from $\Or(n,\C)$ to $\Or(n)$ in the bundle $(V,Q)$ as follows:
choose a Hermitian metric on $V$ and compose the corresponding
isomorphism $V^* \to \bar{V}$ with the isomorphism $q\colon V \to V^*$
to obtain the real structure on $V$.  We shall take a hermitian metric
under which the $F_{i}$ are orthogonal.

When $n=2q$ is even, we see that the underlying real bundle is 
\begin{displaymath}
  F_{1/2} \oplus \cdots \oplus F_{m}
  = F_{-1/2} \otimes (K^{-1} \oplus \cdots \oplus K^{-q})\ .
\end{displaymath}
Recalling that $F_{-1/2}^{2} = K$, this gives the formula stated for 
$w_{2}$.

When $n=2q+1$ is odd, we see that the underlying real bundle is 
\begin{displaymath}
  F_{0,\R} \oplus F_{1} \oplus \cdots \oplus F_{m}\ ,
\end{displaymath}
where $F_{0,\R}$ is the real bundle underlying $F_{0}$ in the real
structure defined by $q \colon F_{0} \to F_{0}^{*}$ and the hermitian
metric.  We have that $w_{1}(F_{i}) = 0$ and $w_{2}(F_{i}) =
c_{1}(F_{i}) = 0 \mod 2$, since $F_{i} = F_{0} \otimes K^{-i}$ has
even degree.  Hence $w_{i}(V,Q) = w_{i}(F_{0,\R})$.  To complete the
proof we only need to note that $w_{2}(F_{0,\R}) = c_{1}(F_{0}) \mod 2
= 0$.
\end{proof}

One of the main results of \cite{hitchin:1992} is the construction of
a ``Teichm\"uller component'' of the moduli space of representations
of $\pi_{1}(X)$ in any split real form of a complex simple Lie group.
The construction in the case of representations in $\SL(n,\R)$ is
quite explicit (see \cite[\S 3]{hitchin:1992}): one keeps the same
underlying bundle and adds certain extra entries to the matrix of the
Higgs field, thus parametrizing the component by spaces of sections of
powers of $K$.  The construction carries over to the case of
$\GL(n,\R)$-Higgs bundles with only one small modification: in
\cite[(3.2)]{hitchin:1992}, substitute the zero in the bottom row of
the matrix with $\alpha_0 \in H^0(X,K)$.  We thus have the following
result.

\goodbreak

\begin{proposition}\label{prop:teich-comp}
  Let $(V,Q,\Phi)$ be a local minimum of $f$ of the kind described in
  $(2)$ of Theorem~\ref{thm:GL-minima} with invariants $(w_1,w_2)$.
  Then there is a connected component of $\mathcal{M}(w_1,w_2)$
  containing $(V,Q,\Phi)$ and isomorphic to a vector space. Any
  $(V,Q,\Phi)$ in this component has $\Phi \neq 0$.
\end{proposition}

%%%%%%%%%%%%%%%%%%%%%%%%%%%%%%%%%%%%%%%%%%%%%%%%%%%%%%%%%%%%%%%%%%%%%%
\section{Connected components of the moduli space}
%%%%%%%%%%%%%%%%%%%%%%%%%%%%%%%%%%%%%%%%%%%%%%%%%%%%%%%%%%%%%%%%%%%%%%

In this section we finally determine the connected components of
$\mathcal{M}^{+}_{\GL(n,\R)}$.
In order to deal with possible non-stable minima we need the following
result.  The proof given at the beginning of \S 10 of
\cite{hitchin:1992} for $G=\SL(n,\R)$ also works for $G = \GL(n,\R)$.

\begin{lemma}\label{lem:unstable-minima}
  Any $(V,Q,\Phi)$ representing a local minimum of $f$ and which has 
  $\Phi \neq 0$ is a stable $\GL(n,\R)$-Higgs bundle. 
\end{lemma}

We can now state and prove our main Theorem.

\begin{theorem}
  \label{thm:main}
  Let $\mathcal{M}^{+}_{\GL(n,\R)}$ be the moduli space of semi-simple
  representations in $\GL(n,\R)$ of the fundamental group of a closed
  oriented surface of genus $g \geq 2$.  Assume that $n \geq 3$. For
  $(w_{1},w_{2}) \in H^{1}(X,\Z/2) \oplus H^{2}(X,\Z/2)$ let
  $\mathcal{M}^{+}_{\GL(n,\R)}(w_{1},w_{2})$ be the subspace of
  representations with invariants $(w_{1},w_{2})$.
  \begin{itemize}
    \item[$(1)$] If $n=2q$ is even, then $\mathcal{M}^{+}_{\GL(n,\R)}$ 
    has $3\cdot 2^{2g}$ connected components.  More precisely:
    \begin{itemize}
      \item[$(i)$]  If $w_{1} \neq 0$, then 
      $\mathcal{M}^{+}_{\GL(n,\R)}(w_{1},w_{2})$ is connected.
      \item[$(ii)$]  If $w_{1} = 0$ and $w_{2} \neq (g-1)q^2 \mod 2$, then 
      $\mathcal{M}^{+}_{\GL(n,\R)}(w_{1},w_{2})$ is connected.
      \item[$(iii)$]  If $w_{1} = 0$ and $w_{2} = (g-1)q^2 \mod 2$, then 
      $\mathcal{M}^{+}_{\GL(n,\R)}(w_{1},w_{2})$ has $2^{2g}+1$ connected 
      components.
    \end{itemize}

    \item[$(2)$] If $n=2q+1$ is odd, then $\mathcal{M}^{+}_{\GL(n,\R)}$ 
    has $3\cdot 2^{2g}$ connected components.  More precisely:
    \begin{itemize}
      \item[$(i)$] If $w_{2} \neq 0$, then
      $\mathcal{M}^{+}_{\GL(n,\R)}(w_{1},w_{2})$ is connected. 
      \item[$(ii)$] If $w_{2} = 0$, then
      $\mathcal{M}^{+}_{\GL(n,\R)}(w_{1},w_{2})$ has $2$ connected
      components.
    \end{itemize}
  \end{itemize}
\end{theorem}

\begin{proof}
By Theorem~\ref{thm:CDHS}, $\mathcal{M}^{+}_{\GL(n,\R)}(w_{1},w_{2})$
has the same number of connected components as
$\mathcal{M}(w_{1},w_{2})$.

Any connected component of this space must contain a minimum of the
non-negative proper map $f$.  Thus, combining
Theorem~\ref{thm:GL-minima} and Lemma~\ref{lem:unstable-minima}, each
component contains either a $\GL(n,\R)$-Higgs bundle $(V,Q,\Phi)$ with
$\Phi = 0$, or one with $\Phi \neq 0$ (i.e.\ of the form given in
$(2)$ of Theorem~\ref{thm:GL-minima}).
Proposition~\ref{prop:teich-comp} shows that for each isomorphism
class of a minimum $(V,Q,\Phi)$ of $f$ with $\Phi \neq 0$, there is a
connected component of the corresponding $\mathcal{M}(w_1,w_2)$.
Since $\Phi \neq 0$ for all $(V,Q,\Phi)$ in this component, it is
disjoint from any component with a minimum with $\Phi = 0$.

For each value of $(w_{1},w_{2})$ we let
$\mathcal{M}_{0}(w_{1},w_{2})$ be the space $\mathcal{M}(w_{1},w_{2})$
with any components with minima with $\Phi \neq 0$ removed.  The space
$\mathcal{M}(w_{1},w_{2})$ contains the moduli space of semistable
principal $\Or(n,\C)$-bundles with the same invariants
$(w_{1},w_{2})$, included as $(V,Q) \mapsto (V,Q,0)$.  From
Theorem~\ref{thm:GL-minima} we conclude that this is exactly the space
of local minima of $f$ on $\mathcal{M}_{0}(w_{1},w_{2})$.  But from
Ramanathan \cite[Proposition~4.2]{ramanathan}, we know that the moduli
space of semistable principal $\Or(n,\C)$-bundles is connected.  Hence
Proposition~\ref{prop:proper} shows that
$\mathcal{M}_{0}(w_{1},w_{2})$ is a connected component.

In conclusion we then have one connected component
$\mathcal{M}_{0}(w_{1},w_{2})$ for each value of $(w_{1},w_{2})$ and
an ``extra'' connected component for each isomorphism
class of a minimum $(V,Q,\Phi)$ of $f$ with $\Phi \neq 0$.  With the aid 
of Propositions~\ref{prop:char-minima} and \ref{prop:inv-minima} this 
leads to the statement of the Theorem.
\end{proof}

\begin{remark}
  When $n$ is even, our count differs somewhat from the count for $G =
  \PSL(n,\R)$ in \cite{hitchin:1992}:
  
  (1) The existence of a different connected component for each choice
  of a square root $F_{-1/2}$ of $K$ (cf.\ $(1)$ of
  Proposition~\ref{prop:char-minima}), giving rise to the $2^{2g}$
  components with $\Phi \neq 0$ in $(iii)$ of $(1)$ of
  Theorem~\ref{thm:main} does not occur in \cite{hitchin:1992}.  This
  is because the Higgs bundles in question are projectively
  equivalent.
  
  (2) It was observed in \cite[p.\ 473]{hitchin:1992} that, when $n$
  is even, the components corresponding to minima with $\Phi \neq 0$
  appear twice in the moduli space of $\PSL(n,\R)$-representations.
  This happens because in that paper $\PSL(n,\R)$-representations
  are analyzed, in the first place, up to $\PGL(n,\R)$-equivalence.
  Since we are dealing here with $\GL(n,\R)$-representations up to
  $\GL(n,\R)$-equi\-va\-lence, this phenomenon does not occur in our case.
\end{remark}

%%%%%%%%%%%%%%%%%%%%%%%%%%%%%%%%%%%%%%%%%%%%%%%%%%%%%%%%%%%%%%%%
\section*{Acknowledgments}
The authors are members of VBAC (Vector Bundles on Algebraic Curves),
which is partially supported by EAGER (EC FP5 Contract no.\ 
HPRN-CT-2000-00099) and by EDGE (EC FP5 Contract no.\ 
HPRN-CT-2000-00101).  The second and third authors were partially
supported by the Portugal/Spain bilateral Programme Acciones
Integradas, grant nos.\ HP2000-0015 and AI-01/24.  The first author
was partially supported by the National Science Foundation under grant
DMS-0072073.  The second author was partially supported by the
Ministerio de Ciencia y Tecnolog\'{\i}a (Spain) under grant
BFM2000-0024. The third author was partially supported by the Funda{\c
  c}{\~a}o para a Ci{\^e}ncia e a Tecnologia (Portugal) through the
Centro de Matem{\'a}tica da Universidade do Porto.

%%%%%%%%%%%%%%%%%%%%%%%%%%%%%%%%%%%%%%%%%%%%%%%%%%%%%%%%%%%%%%%%
%
% Bibliography. Follow the usual conventions. For instance:
%
%%%%%%%%%%%%%%%%%%%%%%%%%%%%%%%%%%%%%%%%%%%%%%%%%%%%%%%%%%%%%%%%

%%%%%%%%%%%%%%%%%%%%%%%%%%%%%%%%%%%%%%%%%%%%%%%%%%%%%%%%%%%%%%%%
\end{document}